\documentclass[11pt]{article}
\usepackage{amsmath,amssymb,amsthm,graphics,latexsym,amsfonts}
\usepackage{fancyhdr}
\usepackage{color}
\usepackage{graphics}

\allowdisplaybreaks[4]

\newtheorem{theorem}{Theorem}
\newtheorem{lemma}[theorem]{Lemma}

\newtheorem{conjecture}[theorem]{Conjecture}

\usepackage{latexsym,amsmath,amssymb,amsfonts,epsfig,graphicx}

\usepackage{color,colordvi,amscd}

\usepackage{verbatim}

\newcommand{\pf}{\noindent {\bf Proof.} }

\begin{document}

\title{A characterization on $(g,f)$-parity orientations\thanks{Supported by National Key R\&D Program of China under grant No. 2023YFA1010203 and  the National Natural
Science Foundation of China under grant No.12271425}}

\author{
 Hongliang Lu and Xinxin Ma\footnote{Corresponding author: luhongliang215@sina.com}\\School of Mathematics and Statistics\\
Xi'an Jiaotong University\\
Xi'an, Shaanxi 710049, China\\
\smallskip\\
}
\date{}

%\underline{}\maketitle

\date{}

\maketitle

\begin{abstract}
Let $G$ be a graph and $g,f:V(G)\to2^N$ be two set functions such that $g(v)\le f(v)\ \mbox{and}\  g(v)\equiv f(v)\pmod 2\ \mbox{for every}\ v\in V(G)$.
An orientation $O$ of $G$ is called a $(g,f)$-parity orientation if $g(v)\le d^+_O(v)\le f(v)$ and $g(v)\equiv  d^+_O(v)\pmod 2$ for every $v\in V(G)$. In this paper, we give a Tutte-type  characterization for a graph to have a $(g,f)$-parity orientation.

\end{abstract}

\section{Introduction}
Let $G$ be a  graph with vertex set $V(G)$ and edge set $E(G)$. Let $e(G):=|E(G)|$. For $v\in V(G)$, we use $d_G(v)$ to denote the degree of a vertex $v$ in $G$. For $S\subseteq V(G)$, let $G[S]$ denote the vertex induced subgraph induced by $S$. For two disjoint vertex sets $S,T$ of $V(G)$, let $E_G(S,T)$ denote the set of edges of $G$ joining $S$ to $T$ and let $e_G(S,T):=|E_G(S,T)|$. When $S=\{x\}$, we denote $E_G(\{x\},T)$ and $e_G(\{x\},T)$ by $E_G(x,T)$ and $e_G(x,T)$, respectively.  Given two  integers $s,t$ such that $s\leq t$, $\{s,s+1,\ldots,t\}$ is denoted by $[s,t]$. An orientation of $G$ is an assignment of a direction to each edge of $G$. For a vertex $v\in V(G)$, we denote by $d^+_O(v)$ the out-degree of $v$ under the orientation $O$ of the edges of $G$. Given an orientation $O$ of $G$ and two disjoint subsets $X,Y$, we use $E_O(X,Y)$  to  denote the set of directed edges of $O$ with tail in $X$ and  head in $Y$ and write $e_O(X,Y):=|E_O(X,Y)|$.

 Because of its fruitful applications, an orientation with specified properties has been extensively studied \cite{NM,JD}. Frank and Gy\'arf\'as \cite{AA} proved that for a graph $G$ and two mappings $a,b: V(G) \rightarrow \mathbb{N}$ with $a(v)\leqslant b(v)$ for every vertex $v$, $G$ has an orientation such that $$a(v)\leqslant d^+(v)\leqslant b(v)$$ for every vertex $v$, if and only if for any $U\subseteq V(G)$, $$\sum_{v\in U}a(v)-d(U)\leqslant |E(G[U])|\leqslant \sum_{v\in U}b(v).$$ where $d(U)$ is the number of edges connecting $U$ and $V(G)\setminus U$. Borowiecki, Grytczuk and Pil\'{s}niak \cite{MJ} showed that every graph
admits an orientation  such that the outdegrees of any two adjacent vertices
are different. These orientations can be interpreted as graph colorings and are now known
as proper orientations. %And the proper orientation number $\vec{\chi}(G)$ of a graph $G$ is the minimum value, taken over all proper orientations of $G$. Chen, Mohar and  Wu \cite{WU} obtained that if $G$ is a planar graph, then $\vec{\chi}(G)\leqslant 14$.

Let $F:V(G)\to2^N$. A graph is said to be \emph{F-avoiding} if there exists an orientation $O$ of $G$ such that $d^+_O(v)\notin F(v)$ for every $v\in V(G)$. In this case, we say that $O$ avoids $F$ or $O$ is an  \emph{F-avoiding} orientation. Conversely, for $H:V(G)\to2^N$, we call an orientation $O$ of $G$ an  \emph{H-orientation}, if $d^+_O(v)\in H(v)$ for every $v\in V(G)$. %In the case, we also say that $G$ admits an $H$-orientation, otherwise, $G$ does not admits  $H$-orientation.
 In this paper,  we may always assume that for every $v\in V(G)$, $\max H(v)\leq d_G(v)$ and  $\min H(v)\geq 0$.

 Akbari et al. \cite{SM} proved that for a graph $G$ and $F:V(G)\to2^N$, if
$$|F(v)|\leq \frac{d_G(v)}{4}$$
for every vertex $v$, then $G$ is \emph{F-avoiding}. Furthermore, they proposed the following conjecture.
\begin{conjecture}[Akbari et al,\cite{SM}]\label{AK}
	For a graph $G$, and a mapping  $F:V(G)\to2^N$, if
	$$|F(v)|\leq \frac{1}{2}(d_G(v)-1)$$
	for every $v\in V(G)$, then $G$ is  \emph{F-avoiding}.
\end{conjecture}
In the same paper, Akbari et al. \cite{SM} proved that Conjecture \ref{AK} holds for bipartite graphs. Bradshaw et.al. made some research progress on Conjecture \ref{AK}: they \cite{BCMMW} proved if $|F(v)|\leq |d_G(v)|/3$ for every $v\in V(G)$, then $G$ is $F$-avoiding. Ma and Lu \cite{ML} gave a characterization for a connected graph $G$ to have an $H$-orientation if $H:V(G)\rightarrow 2^{\mathbb{N}}$ such that
%also made some progress with some constraints on the mapping $F$. %We used the technique of  \emph{parity trace} to describe the existence of the suitable orientation $O$ of $G$ in the constraints of the mapping $H:V(G)\to2^N$, where $H$ has the property:
\begin{align*}
i\notin H(v)\ \mbox{ implies}\  i+1\in H(v),\ \mbox{for every}\  i\in \mathbb{N},\ v\in V(G).
\end{align*}

%Let us call the pair $(G,H)$ \emph{dense} if $H$ satisfies (\ref{main-eq}), and $G$ is connected. \emph{Parity trace} is a specific permutation of all the cut vertices of $G$, the details of which are described in \cite{ML}.
%They gave a Tutte-type characterization for $H$-orientation problem when $(G,H)$ is dense.

%\begin{theorem}\label{main}
%Let $G$ be a graph and let $H:V(G)\rightarrow 2^N$ such that $(G,H)$ is \emph{dense}. There exists an $H$-orientation of $G$ if and only  if there is no parity trace such that $|V^1|\neq e(G)\pmod 2$.
%\end{theorem}

Let $g,f:V(G)\rightarrow \mathbb{N}$ such that $g(v)\equiv f(v)\pmod 2$ for all $v\in V(G)$. An $H$-orientation $O$ of $G$ is called a \emph{$(g,f)$-parity orientation} if $H(v)=\{g(v),g(v)+2,\ldots, f(v)\}$ for all $v\in V(G)$. So a \emph{$(g,f)$-parity orientation $O$} satisfies the following condition: $$g(v)\le d^+_O(v)\le f(v)\ \mbox{and}\ g(v) \equiv f(v)\pmod 2  \mbox{ for every } v\in V(G).$$ A spanning subgraph $F$ of $G$ is called a \emph{$(g,f)$-parity factor} if $$g(v)\le d_F(v)\le f(v)\ \mbox{and}\ d_F(v) \equiv f(v)\pmod 2$$ for every $v\in V(G)$.
 For $U\subseteq V(G)$, we use $f(U),g(U)$ to denote  $\sum_{v\in U}f(v)$ and $\sum_{v\in U}g(v)$, respectively.

% Now in this paper, assume that $H$ satisfies the property:
% \begin{align}\label{main-eqq}
% 	i\notin H(v)\ \mbox{ implies}\  i+1\in H(v),\ \mbox{for}\  g(v)\le i\le f(v).
% \end{align}
%
% where $g(v)=min\{r|r\in H(v)\}$ and $f(v)=max\{r|r\in H(v) \}$.
% In this case, whether $G$ has a \emph{(g,f)-orientation}? To settle this problem, we consider \emph{parity (g,f)-orientations} in which every vertex only has odd outdegree or even outdegree, respectively. These orientations are special cases of $(g,f)$-orientations.
%
% Let $G$ be a general graph, and $g,f:V(G)\rightarrow \mathbb{Z}$ be functions such that $$g(v)\le f(v)\ \mbox{and}\  g(v)\equiv f(v)\pmod 2\ \mbox{for every}\ v\in V(G)$$.
%
% Then an orientation $O$ of $G$ is called a parity $(g,f)$-orientation if $$g(v)\le d^+_O(v)\le f(v)\ \mbox{and}\ g(v) \equiv f(v)\pmod 2$$ for every $v\in V(G)$.  Note that we allow $g(x)\le 0$ and $d_G(y)\le f(y)$ for some
 %In generally, we may assume $0\le g(v)\le f(v)\le d_G(v)$ for all $v\in V(G)$. However, when we apply a criterion for the existence of $(g,f)$-parity orientations, this relaxation on

 In this  paper, we give a characterization for a graph to have a $(g,f)$-parity orientation.
 \begin{theorem}\label{key}
  Let $G$ be a  graph without isolated vertices, and $g,f:V(G)\rightarrow \mathbb{Z}$ be functions such that $g(v)\le f(v)\ \mbox{and}\  g(v)\equiv f(v)\pmod 2\ \mbox{for every}\ v\in V(G)$. Then $G$ has a $(g, f)$-parity orientation if
  and only if for any two disjoint subsets $S$ and $T$ of $V(G)$,
  \begin{align}\label{maineq}
  	\eta_G(S,T)=(f(S)-e_G(S))+(\sum_{v\in T}d_G(v)-e_G(T)-g(T))-q_G(S,T)\ge0
  \end{align}
where $q_G(S,T)$ denotes the number of connected components $C$ of $G-(S\cup T)$ such that
 \begin{align}\label{mainequ}
	f(C)+e_G(C)+e_G(S,C)\equiv 1\pmod 2.
\end{align}
 \end{theorem}

\section{Proof of Theorem \ref{key}}

For completing the proof of Theorem \ref{key}, we need the following   result obtained by Lov\'{a}sz \cite{LA}.
 \begin{theorem}\label{useful}
	Let $G$ be a  graph, and let $g,f:V(G)\rightarrow \mathbb{Z}$ be two functions such that $g(v)\le f(v)\ \mbox{and}\  g(v)\equiv f(v)\pmod 2\ \mbox{for every}\ v\in V(G)$. Then $G$ has a parity factor if
	and only if for all disjoint subsets $S$ and $T$ of $V(G$), it follows that
	\begin{align}\label{tho}
		\delta_G(S,T)=f(S)-g(T)+\sum_{v\in T}d_{G-S}(v)-q_G(S,T;f)\ge0
	\end{align}
	where $q_G(S,T;f)$ denotes the number of components $C$ (also called $f$-odd component) of $G-(S\cup T)$, such that
	\begin{align}
		f(C)+e_G(C,T)\equiv 1\pmod 2
	\end{align}
\end{theorem}

%A component $C$ of $G-(S\cup T)$ is called  an \emph{f-odd component} if it satisfies that $$f(C)+e_G(C,T)\equiv 1\pmod 2.$$
Let $Sub(G)$ denote the graph obtained from $G$  by inserting one new vertex for every  edge of $G$. Let $X:=V(G)$ and let $Y$ denote the set of the new adding vertices. One can see that $Sub(G)$ is a bipartite graph. Note that $d_{Sub(G)}(v)=2$ for every $v\in V(Sub(G))-V(G)$. %and we define $g(v)=f(v)={1}$ for every $v\in Y$. By subdividing $G$, we can obtain the following conclusion.

\begin{lemma} \label{pre}
	 Let $G$ be a graph, and $g,f:V(G)\rightarrow \mathbb{Z}$ be two functions such that $$g(v)\le f(v)\ \mbox{and}\  g(v)\equiv f(v)\pmod 2\ \mbox{for every}\ v\in V(G).$$ Then $G$ has a parity $(g,f)$-orientation if and only if  $Sub(G)$ has a $(g',f')$-parity factor, where
 $g',f':V(Sub(G))\rightarrow \mathbb{Z}$ such that $g'(v)=g(v)$, $f'(v)=f(v)$ for all $v\in V(G)$ and $g'(v)=f'(v)=1$ for all $v\in V(Sub(G))-V(G)$.
\end{lemma}
\pf
Write $X:=V(G)$ and $Y:=Sub(G)-V(G)$. For edge $xy\in E(G)$, we denote the vertex  adjacent with $x$ and $y$ in $Sub(G)$ by $v_{xy}$.

Firstly, we show the necessity. Suppose $G$ has a parity $(g,f)$-orientation $O$. Define
\[
M_O(G):=\{xv_{xy}\ |\ \overrightarrow{xy}\in E(O)\}.
\]
Let $F$ be  a graph with vertex set $X\cup Y$ and edge set $M_O(G)$. Since $O$ is a $(g,f)$-orientation of $G$, $d_F(v)=d_{O^+}(v)\in \{g(v),g(v)+2,\ldots,f(v)\}$ for all $v\in X$. So  for all $v\in X$, $d_F(v)\in \{g(v),g(v)+2,\ldots,f(v)\}$ for all $v\in X$. Note $d_F(v)=1$ for all $v\in V(G)$ by the definition of $M_{O}(G)$. So $F$ is  a $(g',f')$-parity factor of $Sub(G)$.
 %$F(G^\ast)\subseteq E(G^\ast)$ whose directions are from $X$ to $Y$ with respect to $O$ is a parity $(g,f)$-factor of $G^\ast$.

Next we show the sufficiency.
Suppose that $Sub(G)$ has a $(g',f')$-parity factor $F$. Note that $d_F(v)=1$ for all $v\in Y$. So for every edge $xy\in E(G)$, $|E(F)\cap \{xv_{xy},yv_{xy}\}|=1$. Now for any edge $xy\in E(G)$, we  orient $xy$ by $\overrightarrow{xy}$ if $xv_{xy}\in E(F)$, otherwise, the edge $xy$ is oriented by $\overrightarrow{yx}$. This results a orientation  of $G$ denoted by $O$. Since $d_F(v)\in \{g(v),g(v)+2,\ldots,f(v)\}$ for all $v\in X$, then we have $d_F(v)=d_{O^+}(v)\in \{g(v),g(v)+2,\ldots,f(v)\}$. So $O$ is a $(g,f)$-parity orientation.
This completes the proof. \qed
%\begin{equation*}
%  \left\{
%     \begin{array}{ll}
%      \overrightarrow{xy}, & \hbox{$xv_{xy}\in E(F)$;} \\
%       \overrightarrow{yx}, & \hbox{$yv_{xy}\in E(F)$.}
%     \end{array}
%   \right.
%\end{equation*}
%
%
% Then we can obtain an orientation $O$ of $E(G^\ast)$ by orienting the edges of $E(F)$ from $X$ to $Y$, and the edges of $E(G^\ast)\setminus E(F)$ from $Y$ to $X$. It is can be seen that $O$ is a parity $(g,f)$-orientation of $G^\ast$. And we can get a parity $(g,f)$-orientation of $G$ by orienting the edge $uv\in E(G)$ the same as the orientation of the directed 2-path of $G^\ast$ whose start and end point are $u$ and $v$.
\medskip

\noindent \textbf{Proof of Theorem \ref{key}.}
Firstly, we show the necessity. Suppose that $G$ has a $(g,f)$-parity orientation $O$. Let $S,T$ be two disjoint subsets of $V(G)$. Let $C_1,\ldots, C_q$ denote these components $C$ of $G-S-T$ such that $	f(V(C))+e_G(C)+e_G(S,V(C))\equiv 1\pmod 2.$ Note that for $1\leq i\leq q$,
\begin{align*}
f(V(C_i))&\equiv\sum_{v\in V(C_i)}d_{O^+}(v)\ {\pmod 2}\\&
=e_G(C_i)+e_O(V(C_i),S\cup T)\\
&=e_G(C_i)+e_G(S,V(C_i))-e_O(S,V(C_i))+e_O(V(C_i),T),
\end{align*}
which implies that
\begin{align}\label{parity-orien-odd}
e_O(S,V(C_i))+e_O(V(C_i),T)\geq 1.
\end{align}

 For every connected component $C$ of $G-S-T$,  at least one of the following three statements holds:
\begin{itemize}
  \item [(a)] there exist $u\in V(C)$ and $v\in S$ such that $\overrightarrow{vu}\in E(O)$;
  \item [(b)] there exist $u\in V(C)$ and $x\in T$ such that $\overrightarrow{ux}\in E(O)$'
  \item [(c)] $E_O(S,V(C))=\emptyset$ and $E_O(V(C),T)=\emptyset$.
\end{itemize}
For $i\in \{a,b\}$, let $q_i$ denote these  components $C$ of $G-S-T$ such that the statement (i) holds.
By (\ref{parity-orien-odd}), for every  $C_i$, at least one of the statements (a) and (b) holds. Thus we may infer that $q_a+q_b\geq q$.
Since $O$ is a $(g,f)$-parity orientation, then we have
%\begin{subequations}
\begin{align}
  &\sum_{x\in S}d_{O^+}(x)=e_G(S)+e_{O}(S,V(G)-S)\leq f(S), \label{eq_OS}\\ %\tag{eq_OS}\\
  &\sum_{x\in T}d_{O^+}(x)=e_G(T)+e_{O}(T,V(G)-T)\geq g(T). \label{eq_OT}
\end{align}
Combining (\ref{eq_OS}) and (\ref{eq_OT}), we have
\begin{align}\label{eq_OST}
 f(S)-e_G(S)-g(T)+e_G(T)\geq e_{O}(S,V(G)-S)-e_{O}(T,V(G)-T).
\end{align}
%\end{subequations}
On the other hand,
\begin{align}
 e_{O}(S,V(G)-S)&\geq e_O(S,T)+q_a, \label{eq_OS1}\\
  e_G(T,V(G)-T)-q_b-e_O(S,T)&\geq e_{O}(T,V(G)-T).\label{eq_OS2}
\end{align}
Thus we may infer that
\begin{align*}
 e_{O}(S,V(G)-S)&\geq e_O(S,T)+q_a \\
  &=e_G(S,T)- e_O(T,S)+q_a\\
&= e_G(S,T)-e_{O}(T,V(G)-T)+e_{O}(T,V(G)-T-S)+q_a\\
&\geq e_G(S,T)-(e_{G}(T,V(G)-T)-q_b)+e_{O}(T,V(G)-T-S)+e_O(S,T)+q_a\\
&=e_G(S,T)-e_{G}(T,V(G)-T)+(e_{O}(T,V(G)-T-S)+e_O(S,T))+q_a+q_b\\
&\geq e_G(S,T)-e_{G}(T,V(G)-T)+e_{O}(T,V(G)-T)-e_O(T,S)+e_O(S,T)+q\\
&\geq e_{O}(T,V(G)-T)-e_{G}(T,V(G)-T)+q.
\end{align*}
i.e.,
\begin{align}\label{eq_OST2}
 e_{O}(S,V(G)-S)-e_{O}(T,V(G)-T)\geq -e_{G}(T,V(G)-T)+q.
\end{align}
By  (\ref{eq_OST}) and (\ref{eq_OST2}), one can see that
\begin{align*}
&f(S)-e_G(S)-g(T)+e_G(T)\geq -e_{G}(T,V(G)-T)+q,
\end{align*}
i.e.,
\begin{align*}
f(S)-e_G(S)-g(T)+\sum_{v\in T}d_G(v)-e_G(T)-q\geq 0.
\end{align*}
This completes the proof of necessity.

% Suppose that there exist two disjoint subsets $S,T$ such that
%\begin{align}\label{maineq}
%  	\eta(S,T)=(f(S)-e_G(S))+(\sum_{v\in T}d_G(v)-e_G(T)-g(T))-q(S,T)<0
%  \end{align}
%where $q(S,T)$ denotes the number of components $C$ of $G-(S\cup T)$ such that
% \begin{align}\label{mainequ}
%	f(C)+e_G(C)+e_G(C,T)\equiv 1\pmod 2.
%\end{align}

Next we show the sufficiency by contradiction. Suppose that $G$ does not have $(g,f)$-parity orientations. Write $G^*:=Sub(G)$. Let $g_1,f_1:V(G^*)\rightarrow \mathbb{Z}$ be two functions such that $g_1(v)=g(v)$, $f_1(v)=f(v)$ for all $v\in V(G)$ and $g_1(v)=f_1(v)=1$ for all $v\in V(G^*)-V(G)$.
Let $g_2,f_2:V(G^*)\rightarrow \mathbb{Z}$ be two functions such that $g_2(v)=g(v)$, $f_2(v)=f(v)$ for all $v\in V(G)$, $g_2(v)=-1$ and $f_2(v)=3$ for all $v\in V(G^*)-V(G)$. Let $X:=V(G)$ and $Y:=V(G^*)-X$.  By Lemma \ref{pre}, $G^*$ contains no $(g_1,f_1)$-parity factors, which implies that  $G^*$ doesn't contain $(g_2,f_2)$-parity factors.

Thus by Theorem \ref{useful}, there exist two disjoint subsets $S,T\subseteq V(G^*)$  such that $S\cup T\neq \emptyset$ and
\begin{align}\label{eq_f2g2}
	\delta_{G^*}(S,T)=f_2(S)-g_2(T)+\sum_{v\in T}d_{G^*-S}(v)-q_{G^*}(S,T;f_2)<0
\end{align}
where $q_{G^*}(S,T;f_2)$ denotes the number of components $R$ (called $f_2$-odd component) of $G^*-(S\cup T)$ such that
\begin{align*}
	f_2(V(R))+e_{G^*}(V(R),T)\equiv 1\pmod 2,
\end{align*}
Write $q:=q_{G^*}(S,T;f_2)$. We denote these $f_2$-odd components of $G^*-S-T$ by $C_1,\ldots,C_q$.   Let $M:=\cup_{i=1}^q C_i$.% and $M:=\cup_{j=1}^t M_j=G^*-S-T-V(C)$.

Now we choose $S$ and $T$ such that $S\cup T$ is minimal.  %We proceed by proving blow two claims, which give characterizations for $S_0$ and $T_0$.

\medskip
\textbf{Claim 1.}~$T\cap Y= \emptyset$.
\medskip

By contradiction, suppose that $T\cap Y\neq\emptyset$. Let $y\in T\cap Y$. Define $T_1:=T-\{y\}$ and $S_1:=S$.
%\medskip
% \textbf {Case 1.}~$d_{G^*-S-V(D)}(y)\geq 1$.
%%\medskip
% Since $d_{G^*}(y)=2$, we consider the situations where $q_{G^*}(S_1,T_1;f_2)=q-1$ and is the minimum. If $y$ is adjacent to $T$ and $C$, or $y$ is adjacent to $C_m$ and $C_n, m\neq n\in [q] $, it can be concluded that $q_{G^*}(S_1,T_1;f_2)$ is smaller than $q$, and $q_{G^*}(S_1,T_1;f_2)=q-1$. And in other cases it is true that $q_{G^*}(S_1,T_1;f_2)\geq q$. So we have $(q_{G^*}(S_1,T_1;f_2))_{min}=q-1$.
%And it can be concluded that
Then we have
\begin{align*}
  \delta_{G^*}(S_1,T_1)&=f_2(S_1)-g_2(T_1)+\sum_{v\in T_1}d_{G^*-S_1}(v)-q_{G^*}(S_1,T_1;f_2) \\
  &\leq f_2(S)-g_2(T)+g(y)+\sum_{v\in T}d_{G^*-S}(v)-d_{G^*-S}(y)-(q-e_{G^*}(y,V(M)))\\
  &= f_2(S)-g_2(T)+g(y)+\sum_{v\in T}d_{G^*-S}(v)-q+(e_{G^*}(y,V(M))-d_{G^*-S}(y))\\
  &\leq  f_2(S)-g_2(T)-1+\sum_{v\in T}d_{G^*-S}(v)-q\\
&< f_2(S)-g_2(T) +\sum_{v\in T}d_{G^*-S}(v)-q <0
\end{align*}
contradicting to the minimality of $S\cup T$. This completes the proof of Claim 1.

\medskip
\textbf{Claim 2.}~$S\cap Y= \emptyset$.
\medskip

Otherwise, suppose that $S\cap Y\neq\emptyset$. Let $x\in S\cap Y$, $S_2:=S-x$ and $T_2:=T$.
% Since $d_{G^*}(x)=2$, we also consider the situations where $q_{G^*}(S_2,T_2;f_2)=q-1$ and is the minimum. If $x$ is adjacent to $S$ and $C$, or $x$ is adjacent to $C$ and $M$, or $e_G(x,V(C))=2$, it can be concluded that $q_{G^*}(S_2,T_2;f_2)$ is smaller than $q$, and $q_{G^*}(S_2,T_2;f_2)=q-1$. And in other cases it is true that $q_{G^*}(S_2,T_2;f_2)\geq q$. So we have $(q_{G^*}(S_2,T_2;f_2))_{min}=q-1$.
Since $d_{G^*}(x)=2$, we may infer that
\begin{align*}
  \delta_{G^*}(S_2,T_2)&=f_2(S_2)-g_2(T_2)+\sum_{v\in T_2}d_{G^*-S_2}(v)-q_{G^*}(S_2,T_2;f_2) \\
  &\leq f_2(S)-f_2(x)-g_2(T)+\sum_{v\in T}d_{G^*-S}(v)+e_{G^*}(x,T)-(q-e_{G^*}(x,V(M)))\\
  &=f_2(S)-g_2(T)+\sum_{v\in T}d_{G^*-S}(v)-q+(e_{G^*}(x,T)+e_{G^*}(x,V(M))-3)\\
  &\leq f_2(S)-g_2(T)+\sum_{v\in T}d_{G^*-S}(v)-q+(d_{G^*}(x)-3)\\
&< f_2(S)-g_2(T)+\sum_{v\in T}d_{G^*-S}(v)-q<0
,
\end{align*}
a contradiction again. This completes the proof of Claim 2.

For $1\leq i\leq q$, write $D_i=G[X\cap V(C_i)]$ and $D:=\cup_{i=1}^q D_i$. Recall that $G^*$ is a bipartite graph with bipartition $(X,Y)$. By Claims 1 and 2, $S\cup T\subseteq X$. Thus every isolated vertex of $G-S-T$ belongs to $Y$.
Note that if $|V(C_i)|\geq 2$, then $|V(D_i)|\geq 1$.

Let $C_1,\ldots,C_r$ denote these non-trivial $f_2$-odd components of $G^*-S-T$, i.e., for $1\leq i
\leq r$, $|V(C_i)|\geq 2$. Then for $r+1\leq i\leq q$,  $C_i$ is an isolated vertex and $V(C_{i})\subseteq Y$. Note that for every vertex $y\in Y$, if $N_{G^*}(y)\subseteq S$ or  $N_{G^*}(y)\subseteq T$, then $y$ is an isolated vertex of $G^*-S-T$, and satisfies that $f_2(y)+e_{G^*}(T,y)\equiv 1\pmod 2$. Moreover, if $N_{G^*}(y)\cap S\neq \emptyset$ and $N_{G^*}(y)\cap T\neq \emptyset$, then $f_2(y)+e_{G^*}(T,y)\equiv 0\pmod 2$. So we may infer that
\begin{align}\label{q=bound}
q=e_G(S)+e_G(T)+r.
\end{align}

Note that  $Sub(D_i)$ may be obtained from  $C_i$ by deleting some vertices of degree one. So $D_i$ is a connected component  of $G-S-T$.
For $1\leq i\leq r$, we have
\begin{align*}
 f_2(C_i)+e_{G^*}(T,V(C_i))&=f_2(V(D_i))+f_2(V(C_i)\cap Y)+e_G(T,V(D_i))\quad (\mbox{by Claim 1})\\
 &=f_2(V(D_i))+3|V(C_i)\cap Y|+e_G(T,V(D_i))\\
&=f_2(V(D_i))+3(e_G(D_i)+e_G(V(D_i),S\cup T))+e_G(T,V(D_i)\\
&\quad \mbox{(by the definition of $G^*$)}\\
&=f_2(V(D_i))+3(e_G(D_i)+e_G(V(D_i),S))+4e_G(T,V(D_i)\\
&\equiv f(V(D_i))+e_G(D_i)+e_G(V(D_i),S)\equiv 1\pmod 2.
\end{align*}
For $1\leq i\leq r$, since  $f_2(C_i)+e_{G^*}(T,V(C_i))\equiv 1\pmod 2$, we have
\begin{align*}
 f(V(D_i))+e_G(D_i)+e_G(V(D_i),S)\equiv 1\pmod 2.
\end{align*}
Recall that $S\cup T\subseteq V(G)=X$ by Claims 1 and 2. Thus by (\ref{eq_f2g2}), we have
\begin{align*}
	0&>\delta_{G^*}(S,T)=f_2(S)-g_2(T)+\sum_{v\in T}d_{G^*-S}(v)-q_{G^*}(S,T;f_2)\\
&=f(S)-g(T)+\sum_{v\in T}d_{G^*-S}(v)-q_{G^*}(S,T;f_2)\\
&=f(S)-g(T)+\sum_{v\in T}d_{G^*}(v)-q_{G^*}(S,T;f_2)\quad \mbox{(since $S\cup T$ is an independent set of $G^*$)}\\
&=f(S)-g(T)+\sum_{v\in T}d_G(v)-(e_G(S)+e_G(T)+r)\quad \mbox{(by (\ref{q=bound}))}\\
&=\eta_G(S,T),
\end{align*}
a contradiction, where $r$ denotes the number of connected components $D_i$ of $G-S-T$ such that $
f(V(D_i))+e_G(D_i)+e_G(V(D_i))\equiv 1$.
This completes the proof. \qed

\end{document}